\documentclass[amstex,12pt,russian,amssymb]{article}

\usepackage{mathtext}
\usepackage[cp1251]{inputenc}
\usepackage[T2A]{fontenc}
\usepackage[russian]{babel}
\usepackage[dvips]{graphicx}
\usepackage{amsmath}
\usepackage{amssymb}
\usepackage{amsxtra}
\usepackage{latexsym}
\usepackage{ifthen}

\begin{document}

\newcounter{lemma}
\newcommand{\lemma}{\par \refstepcounter{lemma}%
{\bf Лемма \arabic{lemma}.}}

\newcounter{corollary}
\newcommand{\corollary}{\par \refstepcounter{corollary}%
{\bf Следствие \arabic{corollary}.}}

\newcounter{remark}
\newcommand{\remark}{\par \refstepcounter{remark}%
{\bf Замечание \arabic{remark}.}}

\newcounter{theorem}
\newcommand{\theorem}{\par \refstepcounter{theorem}%
{\bf Теорема \arabic{theorem}.}}

\newcounter{proposition}
\newcommand{\proposition}{\par \refstepcounter{proposition}%
{\bf Предложение \arabic{proposition}.}}

\renewcommand{\refname}{\centerline{\bf Список литературы}}

\newcommand{\proof}{{\it Доказательство.\,\,}}

\noindent УДК 517.5

{\bf Е.А.~Севостьянов, С.А.~Скворцов} (Житомирский государственный
университет имени И.~Франко)

{\bf Є.О.~Севостьянов, С.О.~Скворцов} (Житомирський державний
університет імені І.~Франко)

{\bf E.A.~Sevost'yanov, S.A.~Skvortsov} (Zhytomyr Ivan Franko State
University)

\medskip
{\bf О равностепенной непрерывности отображений в случае переменных
областей}

{\bf Про одностайну неперервність відображень у випадку змінних
областей}

{\bf On equicontinuity of mappings in a case of variable domains}

\medskip\medskip
Изучаются вопрос о локальном поведении одного класса отображений в
замыкании области евклидового $n$-мерного пространства в случае,
когда отображённая область не является фиксированной. При
определённых условиях на измеримую функцию, определяющую поведение
указанных отображений, а также ограничениях на отображённые области,
установлена равностепенная непрерывность соответствующего семейства
отображений в замыкании исходной области.

\medskip\medskip
Вивчається питання про локальну поведінку одного класу відображень в
замиканні області евклідового $n$-вимірного простору в випадку, коли
відображена область не є фіксованою. За певних умов на вимірну
функцію, яка визначає поведінку вказаних відображень, а також
обмеженнях на відображені області, встановлено одностайну
неперервність відповідної сім'ї відображень в замиканні вихідної
області.

\medskip\medskip
We study a local behavior of one class of mappings, which are
defined in a domain of $n$-measured Euclidean space, in a case, when
corresponding images of this domain are variable. Under some
conditions on a function defining a behavior of mappings mentioned
above, and some restrictions on mapped domains, the equicontinuity
of the corresponding family of the mappings in the closure of the
initial domain is proved.

\newpage
{\bf 1. Введение.} Основные определения и обозначения, встречающиеся
в тексте, могут быть найдены в монографиях \cite{MRSY} и
\cite{Sev$_7$}.

В относительно недавней работе \cite{Sev$_2$} установлена
равностепенная непрерывность одного семейства пространственных
гомеоморфизмов с неограниченной характеристикой квазиконформности,
переводящих заданную область $D$ на заданную область $D^{\,\prime}.$
Отметим, что в \cite{Sev$_2$} поведение отображений связано с их
непрерывным продолжением на границу в поточечном смысле. Основная
цель настоящей статьи -- распространить те же результаты на случай,
когда отображённые области могут меняться и зависят от
фиксированного отображения данного семейства. Отдельно будет
рассмотрен также случай отображений с ветвлением и случай, когда
непрерывное продолжение на границу и равностепенную непрерывность
следует понимать в терминах простых концов. По этому поводу см.
также классические результаты Р.~Някки и Б.~Палка для
квазиконформных отображений \cite{NP}. Теория простых концов для
пространственных областей и граничное поведение квазиконформных
отображений в их терминах развиты Някки \cite{Na}. В случае
отображений с неограниченной характеристикой этот подход был
использован и развит в работах В.Я.~Гутлянского, В.И.~Рязанова,
Д.А.~Ковтонюка и Э.~Якубова (см. \cite{GRY}--\cite{KR}).

\medskip Здесь и далее
\begin{equation}\label{eq49***}
A(x_0, r_1,r_2): =\left\{ x\,\in\,{\Bbb R}^n:
r_1<|x-x_0|<r_2\right\}\,,
\end{equation}
а $M_p(\Gamma)$ означает $p$-модуль семейства кривых $\Gamma.$
Введём в рассмотрение следующее понятие, см. \cite[разд.~7.6
гл.~7]{MRSY}.
Пусть $p\geqslant 1$ и $Q:{\Bbb R}^n\rightarrow [0, \infty]$ --
измеримая по Лебегу функция, $Q(x)\equiv 0$ при всех $x\not\in D.$
Говорят, что отображение $f:D\rightarrow \overline{{\Bbb R}^n}$ есть
{\it кольцевое $Q$-отоб\-ра\-же\-ние в точке $x_0\in \overline{D}$
относительно $p$-модуля,} $x_0\ne \infty,$ если для некоторого
$r_0=r(x_0)$ и произвольных сферического кольца (\ref{eq49***}) и
любых континуумов $E_1\subset \overline{B(x_0, r_1)}\cap D,$
$E_2\subset \left(\overline{{\Bbb R}^n}\setminus B(x_0,
r_2)\right)\cap D,$ отображение $f$ удовлетворяет соотношению
\begin{equation}\label{eq3*!!}
 M_p\left(f\left(\Gamma\left(E_1,\,E_2,\,D\right)\right)\right)\ \leqslant
\int\limits_{A} Q(x)\cdot \eta^p(|x-x_0|)\ dm(x) \end{equation}
для каждой измеримой функции $\eta : (r_1,r_2)\rightarrow [0,\infty
],$ такой что
\begin{equation}\label{eq28*}
\int\limits_{r_1}^{r_2}\eta(r)\ dr\ \geqslant\ 1\,.
\end{equation}
Аналогично, условимся говорить, что отображение $f:D\rightarrow
\overline{{\Bbb R}^n}$ является кольцевым $Q$-отоб\-ра\-же\-нием в
$\overline{D}$ относительно $p$-модуля, если условие (\ref{eq3*!!})
выполнено в каждой точке $x_0\in \overline{D}.$ В точке $x_0=\infty$
данное определение может быть переформулировано при помощи инверсии:
$\varphi(x)=\frac{x}{|x|^2},$ $\infty\mapsto 0.$

\medskip
В дальнейшем $h(x, y)$ обозначает хордальное расстояние между токами
$x, y\in \overline{{\Bbb R}^n},$ а $h(E)$ -- хордальный диаметр
множества $E\subset\overline{{\Bbb R}^n}$ (см. \cite[гл.~1]{MRSY}).
Пусть $D_i$ -- некоторая фиксированная последовательность областей,
$i=1,2,\ldots, .$ Согласно \cite[разд.~2.4]{NP}, будем говорить, что
семейство областей $D_i,$ $i\in I,$  является {\it равностепенно
равномерной относительно $p$-модуля}, если для каждого $r>0$
существует число $\delta>0$ такое, что неравенство
\begin{equation}\label{eq17***}
M_p(\Gamma(F^{\,*},F, D_i))\geqslant \delta
\end{equation}
выполнено для всех $i\in I$ и произвольных континуумов $F,
F^*\subset D,$ подчинённых условиям $h(F)\geqslant r$ и
$h(F^{\,*})\geqslant r.$ Отметим, что в работе \cite{NP}
рассматривался частный случай $p=n,$ а неравенство (\ref{eq17***})
здесь требовалось не только для континуумов $F$ и $F^{\,*},$ но и
произвольных связных множеств. Указанное обстоятельство
(несущественно) отличает приведенное выше определение от
аналогичного определения Някки и Палка (\cite{NP}). Заметим также,
что для фиксированной области $D_i$ соотношение (\ref{eq17***})
влечёт так называемую сильную достижимость её границы относительно
$p$-модуля (см. \cite[теорема~6.2]{Na73}).

\medskip
Для $p\geqslant 1,$ заданного числа $\delta\geqslant 0,$
фиксированной области $D\subset {\Bbb R}^n,$ $n\geqslant 2,$
континуума $A\subset D$ и заданной функции $Q:D\rightarrow[0,
\infty]$ обозначим через $\frak{F}_{Q, A, p, \delta}(D)$ семейство
всех кольцевых $Q$-гомеоморфизмов $f:D\rightarrow \overline{{\Bbb
R}^n}$ в $\overline{D}$ относительно $p$-модуля, удовлетворяющих
условиям $h(f(A))\geqslant\delta$ и $h(\overline{{\Bbb
R}^n}\setminus f(D))\geqslant \delta.$ Полагаем
$$q_{x_0}(r):=\frac{1}{\omega_{n-1}r^{n-1}}\int\limits_{|x-x_0|=r}Q(x)\,dS\,,$$
где $dS$ -- элемент площади поверхности $S,$ и
$$q^{\,\prime}_{b}(r):=\frac{1}{\omega_{n-1}r^{n-1}}\int\limits_{|x-b|=r}Q^{\,\prime}(x)\,dS\,,$$
$Q^{\,\prime}(x)=\max\{Q(x), 1\}.$ Имеют место следующие
утверждения.

\medskip
\begin{theorem}\label{th3} {\sl\, Предположим, $p\in (n-1, n],$ область $D$ локально связна в каждой точке
$x_0\in\partial D,$ и что области $D_f^{\,\prime}=f(D)$ являются
равностепенно равномерными относительно $p$-модуля по всем
$f\in\frak{F}_{Q, A, p, \delta}(D).$ Пусть при $p=n$ число
$\delta>0,$ а при $n-1<p<n$ число $\delta\geqslant 0.$ Если функция
$Q$ имеет конечное среднее колебание в $\overline{D},$ либо в каждой
точке $x_0\in \overline{D}$ при некотором $\beta(x_0)>0$ выполнено
условие
\begin{equation}\label{eq2}
\int\limits_{0}^{\beta(x_0)}\frac{dt}{t^{\frac{n-1}{p-1}}q_{x_0}^{\,\prime\,\frac{1}{p-1}}(t)}=\infty\,,
\end{equation}
то каждое из отображений $f\in\frak{F}_{Q, A, p, \delta}(D)$ имеет
непрерывное продолжение в $\overline{D}$ и семейство $\frak{F}_{Q,
A, p, \delta}(\overline{D}),$ состоящее из всех, таким образом,
продолженных отображений $\overline{f}: \overline{D}\rightarrow
\overline{{\Bbb R}^n},$ является равностепенно непрерывным в
$\overline{D}.$
  }
\end{theorem}

\medskip
В случае отображений с ветвлением теорема \ref{th3} допускает
некоторое обобщение, в связи с чем напомним следующее определение.
Отображение $f:D\rightarrow\overline{{\Bbb R}^n}$ области
$D\subset{\Bbb R}^n$ на область $D^{\,\prime}\subset\overline{{\Bbb
R}^n}$ назовём {\it замкнутым}, если $C(f, \partial D)\subset
\partial D^{\,\prime},$ где, как обычно, $C(f, \partial D)$ --
предельное множество отображения $f$ на $\partial D.$ Для
$p\geqslant 1,$ фиксированной области $D\subset {\Bbb R}^n,$
множества $E\subset\overline{{\Bbb R}^n}$ и числа $\delta>0$
обозначим через $\frak{R}_{Q, \delta, p, E}(D)$ семейство всех
открытых дискретных замкнутых кольцевых $Q$-отображений
$f:D\rightarrow \overline{{\Bbb R}^n}\setminus E$ относительно
$p$-модуля в $\overline{D}$ со следующим условием: для всякого $f$ и
всякой области $D^{\,\prime}_f:=f(D)$ найдётся континуум $K_f\subset
D^{\,\prime}_f$ такой, что $h(K_f)\geqslant \delta$ и
$h(f^{\,-1}(K_f),
\partial D)\geqslant \delta>0.$

\medskip
\begin{theorem}\label{th4} {\sl\, Предположим, $p\in (n-1, n],$ область $D$ локально связна в каждой точке
$x_0\in\partial D$ и что области $D_f^{\,\prime}=f(D)$ являются
равностепенно равномерными относительно $p$-модуля по всем
$\frak{R}_{Q, \delta, p, E}(D).$ Пусть при $p=n$ множество $E$ имеет
положительную ёмкость, а при $n-1<p<n$ является произвольным
замкнутым множеством. Если функция $Q$ имеет конечное среднее
колебание в $\overline{D},$ либо в каждой точке $x_0\in
\overline{D}$ при некотором $\beta(x_0)>0$ выполнено условие
(\ref{eq2}), то каждое из отображений $\frak{R}_{Q, \delta, p,
E}(D)$ имеет непрерывное продолжение в $\overline{D}$ и семейство
$\frak{R}_{Q, \delta, p, E}(\overline{D}),$ состоящее из всех, таким
образом, продолженных отображений $\overline{f}:
\overline{D}\rightarrow \overline{{\Bbb R}^n},$ является
равностепенно непрерывным в $\overline{D}.$ }
\end{theorem}

\medskip
Теоремы \ref{th3} и \ref{th4} относятся к случаю локально связных
границ, когда граничное продолжение отображений и равностепенную
непрерывность их семейств необходимо понимать в обычном,
<<поточечном>> смысле. Случай более сложных границ соответствуют
ситуации так называемых простых концов (определение и
соответствующая терминология может быть найдены в работе \cite{KR}).
Будем говорить, что граница области $D$ в ${\Bbb R}^n$ является {\it
локально квазиконформной}, если каждая точка $x_0\in\partial D$
имеет окрестность $U$, которая может быть отображена квазиконформным
отображением $\varphi$ на единичный шар ${\Bbb B}^n\subset{\Bbb
R}^n$ так, что $\varphi(\partial D\cap U)$ является пересечением
${\Bbb B}^n$ с координатной гиперплоскостью. Говорим, что
ограниченная область $D$ в ${\Bbb R}^n$ {\it регулярна}, если $D$
может быть квазиконформно отображена на область с локально
квазиконформной границей. Если $\overline{D}_P$ является пополнением
регулярной области $D$ ее простыми концами и $g_0$ является
квазиконформным отображением области $D_0$ с локально
квазиконформной границей на $D$, то оно естественным образом
определяет в $\overline{D}_P$ некоторую метрику (см. \cite{KR}).
Справедливо следующее утверждения.

\medskip
\begin{theorem}\label{th1} {\sl\, Предположим, $p\in (n-1, n],$ область $D$ регулярна и что области
$D_f^{\,\prime}=f(D)$ являются ограниченными равностепенно
равномерными относительно $p$-модуля по $f\in\frak{F}_{Q, A, p,
\delta}(D)$ областями с локально квазиконформной границей. Пусть при
$p=n$ число $\delta>0,$ а при $n-1<p<n$ число $\delta\geqslant 0.$
Если функция $Q$ имеет конечное среднее колебание в $\overline{D},$
либо в каждой точке $x_0\in \overline{D}$ при некотором
$\beta(x_0)>0$ выполнено условие (\ref{eq2}), то каждое из
отображений $f\in\frak{F}_{Q, A, p, \delta}(D)$ имеет непрерывное
продолжение $\overline{f}: \overline{D}_P\rightarrow \overline{{\Bbb
R}^n}$ в $\overline{D}_P$ и семейство $\frak{F}_{Q, A, p,
\delta}(\overline{D}),$ состоящее из всех, таким образом,
продолженных отображений $\overline{f}: \overline{D}_P\rightarrow
\overline{{\Bbb R}^n},$ является равностепенно непрерывным в
$\overline{D}_P.$
  }
\end{theorem}

\medskip
\begin{theorem}\label{th2} {\sl\, Предположим, $p\in (n-1, n],$ область $D$ регулярна и что области
$D_f^{\,\prime}=f(D)$ являются ограниченными равностепенно
равномерными относительно $p$-модуля по $\frak{R}_{Q, \delta, p,
E}(D)$ областями с локально квазиконформной границей. Пусть при
$p=n$ множество $E$ имеет положительную ёмкость, а при $n-1<p<n$
является произвольным замкнутым множеством. Если функция $Q$ имеет
конечное среднее колебание в $\overline{D},$ либо в каждой точке
$x_0\in \overline{D}$ при некотором $\beta(x_0)>0$ выполнено условие
(\ref{eq2}), то каждое из отображений $\frak{R}_{Q, \delta, p,
E}(D)$ имеет непрерывное продолжение в $\overline{D}$ и семейство
$\frak{R}_{Q, \delta, p, E}(\overline{D}),$ состоящее из всех, таким
образом, продолженных отображений $\overline{f}:
\overline{D}_P\rightarrow \overline{{\Bbb R}^n},$ является
равностепенно непрерывным в $\overline{D}_P.$ }
\end{theorem}

\medskip
{\bf 2. Формулировка и доказательство основных лемм.} Основным
инструментом на пути доказательства теорем \ref{th3} и \ref{th4}
являются следующие два утверждения.

\medskip
\begin{lemma}\label{lem1}{\sl\, Предположим, область $D$ локально связна в каждой точке
$x_0\in\partial D,$ и что области $D_f^{\,\prime}=f(D)$ являются
равностепенно равномерными относительно $p$-модуля по всем
$f\in\frak{F}_{Q, A, p, \delta}(D).$ Пусть при $p=n$ число
$\delta>0,$ а при $n-1<p<n$ число $\delta\geqslant 0.$ Предположим
также, что для каждой точки $x_0\in \overline{D}$ найдётся
$\varepsilon_0=\varepsilon_0(x_0)>0$ и измеримая по Лебегу функция
$\psi(t):(0, \varepsilon_0)\rightarrow [0,\infty]$ со следующим
свойством: для любого $\varepsilon\in(0, \varepsilon_0)$ выполнено
условие
\begin{equation}\label{eq7***} I(\varepsilon,
\varepsilon_0):=\int\limits_{\varepsilon}^{\varepsilon_0}\psi(t)dt <
\infty\,,\quad I(\varepsilon, \varepsilon_0)\rightarrow
\infty\quad\text{при}\quad\varepsilon\rightarrow 0\,,
\end{equation}
и, кроме того, при  $\varepsilon\rightarrow 0$
\begin{equation} \label{eq3.7.2}
\int\limits_{A(x_0, \varepsilon, \varepsilon_0)}
Q(x)\cdot\psi^{\,p}(|x-x_0|)\,dm(x) = o(I^p(\varepsilon,
\varepsilon_0))\,,\end{equation}
где, как обычно, сферическое кольцо $A(x_0, \varepsilon,
\varepsilon_0)$ определено как в (\ref{eq49***}).
Тогда каждое из отображений $f\in\frak{F}_{Q, A, p, \delta}(D)$
имеет непрерывное продолжение в $\overline{D}$ и семейство
$\frak{F}_{Q, A, p, \delta}(\overline{D}),$ состоящее из всех, таким
образом, продолженных отображений $\overline{f}:
\overline{D}\rightarrow \overline{{\Bbb R}^n},$ является
равностепенно непрерывным в $\overline{D}.$
  }
\end{lemma}

\medskip
\begin{proof}
Равностепенная непрерывность внутри области $D$ вытекает из
\cite[лемма~3.2.2]{Sev$_7$} в случае $p=n$ и \cite[лемма~2.4]{GSS}
при $n-1<p<n,$ а возможность продолжения каждого элемента $f$
семейства отображений $\frak{F}_{Q, A, p, \delta}(D)$ до
непрерывного отображения в замыкании $D$ ~--- из
\cite[лемма~1]{Sev$_4$} при $p=n$ (доказательство этого факта в
случае $n-1<p<n$ проводится аналогично).

Осталось показать, что семейство $\frak{F}_{Q, A, p,
\delta}(\overline{D})$ равностепенно не\-прерывно в точках $\partial
D.$ Предположим противное, тогда найдётся $x_0\in
\partial D$ и число $a>0$ такое, что для каждого
$m=1,2,\ldots$ существуют точка $x_m\in \overline{D}$ и элемент
${f}_m$ семейства $\frak{F}_{Q, A, p, \delta}(\overline{D})$ такие,
что $|x_0-x_m|< 1/m$ и $h(f_m(x_m), f_m(x_0))\geqslant a.$ Поскольку
$f_m$ имеет непрерывное продолжение в точку $x_0,$ то мы можем найти
такую последовательность $x^{\,\prime}_m\in D,$
$x^{\,\prime}_m\rightarrow x_0$ при $m\rightarrow\infty$ такую, что
$h(f_m(x^{\,\prime}_m), f_m(x_0))\leqslant 1/m.$ Таким образом,
\begin{equation}\label{eq6***}
h(f_m(x_m), f_m(x^{\,\prime}_m))\geqslant a/2\qquad \forall\,\,m\in
{\Bbb N}\,.
\end{equation}
%
Можно считать, что $x_0\ne \infty.$ В виду возможности непрерывного
продолжения каждого $f_m$ на границу $D,$ мы можем считать, что
$x_m\in D.$

В силу локальной связности области $D$ в точке $x_0$ найдётся
последовательность окрестностей $V_m$ точки $x_0$ с ${\rm
diam}\,V_m\rightarrow 0$ при $m\rightarrow\infty,$ такие что
множества $D\cap V_m$ являются областями и $D\cap V_m \subset B(x_0,
2^{\,-m}).$ Не ограничивая общности рассуждений, переходя к
подпоследовательности, если это необходимо, мы можем считать, что
$x_m, x^{\,\prime}_m \in D\cap V_m.$ Соединим точки $x_m$ и
$x^{\,\prime}_m$ непрерывной кривой $\gamma_m(t):[0,1]\rightarrow
{\Bbb R}^n$ такой, что $\gamma_m(0)=x_m,$
$\gamma_m(1)=x^{\,\prime}_m$ и $\gamma_m(t)\in V_m$ при $t\in
(0,1).$ Обозначим через $C_m$ образ кривой $\gamma_m(t)$ при
отображении $f_m.$ Из соотношения (\ref{eq6***}) вытекает, что
\begin{equation}\label{eq5.1}
h(C_m)\geqslant a/2\qquad\forall\, m\in {\Bbb N}\,,
\end{equation}
где $h$ обозначает хордальный диаметр множества.

\medskip Не ограничивая общности рассуждений, можно
считать, что континуум $A,$ участвующий в определении класса
$\frak{F}_{Q, A, p, \delta}(D),$ лежит вне шаров $B(x_0, 2^{\,-m}),$
$m=1,2,\ldots, .$ Более того, не ограничивая общности рассуждений,
мы можем считать, что $B(x_0, \varepsilon_0)\cap A=\varnothing.$
Пусть $\Gamma_m$ -- семейство кривых, соединяющих $\gamma_m$ и $A$ в
$D.$ Из определения кольцевого $Q$-отображения относительно
$p$-модуля в точке $x_0$ вытекает, что
\begin{equation}\label{eq10}
M_p(f_m(\Gamma_m))\leqslant \int\limits_{A(x_0, \frac{1}{2^m},
\varepsilon_0)} Q(x)\cdot \eta^p(|x-x_0|)\ dm(x)
\end{equation}
для каждой измеримой функции $\eta:(\frac{1}{2^m},
\varepsilon_0)\rightarrow [0,\infty ],$ такой что
$\int\limits_{\frac{1}{2^m}}^{\varepsilon_0}\eta(r)dr \geqslant 1.$
Заметим, что функция
$$\eta(t)=\left\{
\begin{array}{rr}
\psi(t)/I(2^{\,-m}, \varepsilon_0), &   t\in (2^{\,-m},
\varepsilon_0),\\
0,  &  t\in {\Bbb R}\setminus (2^{\,-m}, \varepsilon_0)\,,
\end{array}
\right. $$ где $I(\varepsilon,
\varepsilon_0):=\int\limits_{\varepsilon}^{\varepsilon_0}\psi(t)dt,$
удовлетворяет условию нормировки вида (\ref{eq28*}) при
$r_1:=2^{\,-m},$ $r_2:=\varepsilon_0,$ поэтому из условий
(\ref{eq3.7.2}) и (\ref{eq10}) вытекает, что
\begin{equation}\label{eq11}
M_p(f_m(\Gamma_m))\leqslant \alpha(2^{\,-m})\rightarrow 0
\end{equation}
при $m\rightarrow \infty,$ где $\alpha(\varepsilon)$ -- некоторая
неотрицательная функция, стремящаяся к нулю при
$\varepsilon\rightarrow 0,$ которая существует ввиду условия
(\ref{eq3.7.2}).

\medskip
С другой стороны, заметим, что $f_m(\Gamma_m)=\Gamma(C_m, f_m(A),
D_m^{\,\prime}).$ По условию леммы $h(f_m(A))\geqslant \delta$ при
всех $m\in {\Bbb N}.$ Следовательно, ввиду (\ref{eq5.1})
$h(f_m(A))\geqslant \delta_1$ и $h(C_m)\geqslant\delta_1,$ где
$\delta_2:=\min\{\delta, a/2\}.$ Воспользовавшись тем, что области
$D_m^{\,\prime}:=f_m(D)$ являются равностепенно равномерными
относительно $p$-модуля, мы заключаем, что существует $\sigma>0$
такое, что
$$M_p(f_m(\Gamma_m))=M_p(\Gamma(C_m, f_m(A),
D_m^{\,\prime}))\geqslant \sigma\qquad\forall\,\, m\in {\Bbb N}\,,$$
что противоречит условию (\ref{eq11}). Полученное противоречие
указывает на то, что предположение об отсутствии равностепенной
непрерывности семейства $\frak{F}_{Q, A, p, \delta}(\overline{D})$
было неверным. Полученное противоречие завершает доказательство
леммы.~$\Box$
\end{proof}

\medskip
В случае отображений с ветвлением лемма \ref{lem1} принимает
следующий вид.

\medskip
\begin{lemma}\label{lem2}{\sl\, Предположим, $p\in (n-1, n],$ область $D$ локально связна в каждой точке
$x_0\in\partial D$ и что области $D_f^{\,\prime}=f(D)$ являются
равностепенно равномерными относительно $p$-модуля по всем
$\frak{R}_{Q, \delta, p, E}(D).$ Пусть при $p=n$ множество $E$ имеет
положительную ёмкость, а при $n-1<p<n$ является произвольным
замкнутым множеством. Предположим также, что для каждой точки
$x_0\in \overline{D}$ найдётся $\varepsilon_0=\varepsilon_0(x_0)>0$
и измеримая по Лебегу функция $\psi(t):(0, \varepsilon_0)\rightarrow
[0,\infty]$ со следующим свойством: для любого $\varepsilon\in(0,
\varepsilon_0)$ выполнено условие (\ref{eq7***}) и, кроме того, при
$\varepsilon\rightarrow 0$ выполнено условие (\ref{eq3.7.2}). Тогда
каждое из отображений $\frak{R}_{Q, \delta, p, E}(D)$ имеет
непрерывное продолжение в $\overline{D}$ и семейство $\frak{R}_{Q,
\delta, p, E}(\overline{D}),$ состоящее из всех, таким образом,
продолженных отображений $\overline{f}: \overline{D}\rightarrow
\overline{{\Bbb R}^n},$ является равностепенно непрерывным в
$\overline{D}.$
  }
\end{lemma}

\medskip
\begin{proof}
Равностепенная непрерывность внутри области $D$ вытекает из
\cite[лемма~3.6.1]{Sev$_7$} в случае $p=n$ и \cite[лемма~2.4]{GSS}
при $n-1<p<n,$ а возможность продолжения каждого элемента $f$
семейства отображений $\frak{R}_{Q, \delta, p, E}(D)$ до
непрерывного отображения в замыкании $D$ ~--- из
\cite[лемма~1]{Sev$_4$} при $p=n$ (доказательство этого факта в
случае $n-1<p<n$ проводится аналогично).

Осталось показать, что семейство $\frak{R}_{Q, \delta, p,
E}(\overline{D})$ равностепенно не\-прерывно в точках $\partial D.$
Предположим противное, тогда найдётся $x_0\in
\partial D$ и число $a>0$ такое, что для каждого
$m=1,2,\ldots$ существуют точка $x_m\in \overline{D}$ и элемент
${f}_m$ семейства $\frak{R}_{Q, \delta, p, E}(\overline{D})$ такие,
что $|x_0-x_m|< 1/m$ и выполнено условие (\ref{eq6***}). Можно
считать, что $x_0\ne \infty.$ В виду возможности непрерывного
продолжения каждого $f_m$ на границу $D,$ мы можем считать, что
$x_m\in D.$ Более того, ввиду того, что $f_m$ продолжается по
непрерывности в точку $x_0,$ найдётся последовательность
$x^{\,\prime}_m\in D,$ сходящаяся к точке $x_0$ при $m\rightarrow
\infty,$ такая, что при некотором $a>0$ выполнены неравенства в
(\ref{eq6***}). Соединим точки $x_m$ и $x^{\,\prime}_m$ кривой
$\gamma_m:[0,1]\rightarrow {\Bbb R}^n$ такой, что $\gamma_m(0)=x_m,$
$\gamma_m(1)=x^{\,\prime}_m$ и $\gamma_m\in V_m\cap D$ при $t\in
(0,1).$ Обозначим через $C_m$ образ кривой $\gamma_m$ при
отображении $f_m.$ Из соотношения (\ref{eq6***}) вытекает, что
выполнено условие вида (\ref{eq5.1}), где $h$ обозначает хордальный
диаметр множества.

По определению семейства отображений $\frak{R}_{Q, \delta, p, E}(D)$
для всякого $f_m$ и всякой области $D^{\,\prime}_m:=f_m(D)$ найдётся
континуум $K_m\subset D^{\,\prime}_m$ такой, что $h(K_m)\geqslant
\delta$ и $h(f_m^{\,-1}(K_m),
\partial D)\geqslant \delta>0.$ Поскольку по условию леммы области $D^{\,\prime}_m$ являются
равностепенно равномерными относительно $p$-модуля, то ввиду
сказанного и учитывая условие (\ref{eq5.1}) мы получим, что при всех
$m=1,2,\ldots$ и некотором $b>0$ выполняется неравенство
\begin{equation}\label{eq13}
M_p(\Gamma(K_m, C_m, D^{\,\prime}_m))\geqslant b\,.
\end{equation}
Рассмотрим семейство $\Gamma_m,$ состоящее из всех кривых $\beta:[0,
1)\rightarrow D^{\,\prime}_m,$ где $\beta(0)\in C_m$ и
$\beta(t)\rightarrow p\in C_m$ при $t\rightarrow 1.$ Пусть
$\Gamma^*_m$ -- семейство всех полных поднятий $\alpha:[0,
1)\rightarrow D$ семейства $\Gamma_m$ при отображении $f_m$ с
началом на $\gamma_m.$ Такое семейство корректно определено ввиду
\cite[теорема~3.7]{Vu}. Ввиду замкнутости отображения $f_m$ имеем:
$\alpha(t)\rightarrow f^{\,-1}_m(K_m),$ где $f^{\,-1}_m(K_m)$ --
полный прообраз континуума $K_m$ при отображении $f_m.$ Не
ограничивая общности рассуждений, можно считать, что $\gamma_m\in
B(x_0, 2^{-m}).$

Заметим, что ввиду компактности пространства $\overline{{\Bbb R}^n}$
при каждом фиксированном $\delta>0$ множество
$C_{\delta}:=\{x\in D: h(x, \partial D)\geqslant \delta\}$
является компактом в $D$ и $f_m^{\,-1}(K_m)\subset C_{\delta}.$
Ввиду \cite[лемма~1]{Smol} множество $C_{\delta}$ можно вложить в
континуум $E_{\delta},$ лежащий в области $D,$ при этом, можно
считать, что ${\rm dist}\,(x_0, E_{\delta})\geqslant \varepsilon_0$
за счёт уменьшения $\varepsilon_0,$ если это необходимо. Тогда на
основании (\ref{eq3*!!}) вытекает, что
\begin{equation}\label{eq10A}
M_p(f_m(\Gamma_m^*))\leqslant M_p(f_m(\Gamma(\gamma_m, E_{\delta},
D)))\leqslant \int\limits_{A(x_0, \frac{1}{2^m}, \varepsilon_0)}
Q(x)\cdot \eta^p(|x-x_0|)\ dm(x)
\end{equation}
для каждой измеримой функции $\eta: (\frac{1}{2^m},
\varepsilon_0)\rightarrow [0,\infty ],$ такой что
$\int\limits_{\frac{1}{2^m}}^{\varepsilon_0}\eta(r)dr \geqslant 1.$
Заметим, что функция
$$\eta(t)=\left\{
\begin{array}{rr}
\psi(t)/I(2^{\,-m}, \varepsilon_0), &   t\in (2^{\,-m},
\varepsilon_0),\\
0,  &  t\in {\Bbb R}\setminus (2^{\,-m}, \varepsilon_0)\,,
\end{array}
\right. $$ где $I(\varepsilon,
\varepsilon_0):=\int\limits_{\varepsilon}^{\varepsilon_0}\psi(t)dt,$
удовлетворяет условию нормировки вида (\ref{eq28*}) при
$r_1:=2^{\,-m},$ $r_2:=\varepsilon_0,$ поэтому из условий
(\ref{eq3.7.2}) и (\ref{eq10A}) вытекает, что
\begin{equation}\label{eq11B}
M_p(f_m(\Gamma_m))\leqslant \alpha(2^{\,-m})\rightarrow 0
\end{equation}
при $m\rightarrow \infty,$ где $\alpha(\varepsilon)$ -- некоторая
неотрицательная функция, стремящаяся к нулю при
$\varepsilon\rightarrow 0,$ которая существует ввиду условия
(\ref{eq3.7.2}). Заметим, кроме того, что $f_m(\Gamma_m)> \Gamma_m$
и, одновременно, $f_m(\Gamma_m)\subset\Gamma_m,$ так что ввиду
\cite[теоремы~6.2, 6.4]{Va}
\begin{equation}\label{eq12}
M_p(f_m(\Gamma_m))=M_p(\Gamma(K_m, C_m, D^{\,\prime}_m))\,.
\end{equation}
Однако, соотношения (\ref{eq11B}) и (\ref{eq12}) в совокупности
противоречат (\ref{eq13}). Полученное противоречие указывает на то,
что исходное предположение (\ref{eq6***}) было неверным, и, значит,
семейство отображений $\frak{R}_{Q, \delta, p, E}(\overline{D})$
равностепенно непрерывно в каждой точке $x_0\in \partial D.$~$\Box$
\end{proof}

\medskip
В случае не локально связных границ заданной области имеют место
следующие аналоги лемм \ref{lem1} и \ref{lem2}.

\medskip
\begin{lemma}\label{lem3} {\sl\, Предположим, $p\in (n-1, n],$ область $D$ регулярна и что области $D_f^{\,\prime}=f(D)$ являются
ограниченными равностепенно равномерными относительно $p$-модуля по
$f\in\frak{F}_{Q, A, p, \delta}(D)$ областями с локально
квазиконформной границей. Пусть при $p=n$ число $\delta>0,$ а при
$n-1<p<n$ число $\delta\geqslant 0.$ Предположим также, что для
каждой точки $x_0\in \overline{D}$ найдётся
$\varepsilon_0=\varepsilon_0(x_0)>0$ и измеримая по Лебегу функция
$\psi(t):(0, \varepsilon_0)\rightarrow [0,\infty]$ со следующим
свойством: для любого $\varepsilon\in(0, \varepsilon_0)$ выполнено
условие (\ref{eq7***}) и, кроме того, при $\varepsilon\rightarrow 0$
выполнено условие (\ref{eq3.7.2}). Тогда каждое из отображений
$f\in\frak{F}_{Q, A, p, \delta}(D)$ имеет непрерывное продолжение
$\overline{f}: \overline{D}_P\rightarrow \overline{{\Bbb R}^n}$ в
$\overline{D}_P$ и семейство $\frak{F}_{Q, A, p,
\delta}(\overline{D}),$ состоящее из всех, таким образом,
продолженных отображений $\overline{f}: \overline{D}_P\rightarrow
\overline{{\Bbb R}^n},$ является равностепенно непрерывным в
$\overline{D}_P.$
  }
\end{lemma}

\begin{proof} Равностепенная непрерывность внутри области $D$ вытекает из
\cite[лемма~3.2.2]{Sev$_7$} в случае $p=n$ и \cite[лемма~2.4]{GSS}
при $n-1<p<n,$ а возможность продолжения каждого элемента $f$
семейства отображений $\frak{F}_{Q, A, p, \delta}(D)$ до
непрерывного отображения в замыкании $D$ ~--- из
\cite[лемма~3]{Sev$_5$}.

Покажем равностепенную непрерывность семейства $\frak{F}_{Q, A, p,
\delta}(\overline{D})$ в точках $E_D,$ где $E_D$ обозначает
пространство простых концов, соответствующее области $D.$
Предположим противное, а именно, что семейство $\frak{F}_{Q, A, p,
\delta}(\overline{D})$ не является равностепенно непрерывным в
некоторой точке $P_0\in E_D.$ Тогда найдутся число $a>0,$
последовательность $P_k\in \overline{D}_P,$ $k=1,2,\ldots$ и
элементы $f_k\in\frak{F}_{Q, A, p, \delta}(D)$ такие, что $d(P_k,
P_0)<1/k$ и
\begin{equation}\label{eq3C}
h(f_k(P_k), f_k(P_0))\geqslant a\quad\forall\quad k=1,2,\ldots,\,.
\end{equation}
Ввиду возможности непрерывного продолжения каждого $f_k$ на границу
$D$ в терминах простых концов, для всякого $k\in {\Bbb N}$ найдётся
элемент $x_k\in D$ такой, что $d(x_k, P_k)<1/k$ и $h(f_k(x_k),
f_k(P_k))<1/k.$ Тогда из (\ref{eq3C}) вытекает, что
\begin{equation}\label{eq4C}
h(f_k(x_k), f_k(P_0))\geqslant a/2\quad\forall\quad k=1,2,\ldots,\,.
\end{equation}
Аналогично, в силу непрерывного продолжения отображения $f_k$ в
$\overline{D_P}$ найдётся последовательность $x_k^{\,\prime}\in D,$
$x_k^{\,\prime}\rightarrow P_0$ при $k\rightarrow \infty$ такая, что
$|f_k(x_k^{\,\prime})-f_k(P_0)|<1/k$ при $k=1,2,\ldots\,.$ Тогда из
(\ref{eq4C}) вытекает, что
\begin{equation}\label{eq5E}
h(f_k(x_k), f_k(x_k^{\,\prime}))\geqslant a/4\quad\forall\quad
k=1,2,\ldots\,,
\end{equation}
где последовательности $x_k$ и $x_k^{\,\prime}$ принадлежат $D$ и
сходятся к простому концу $P_0$ при $k\rightarrow\infty.$

В силу \cite[лемма~2]{KR} простой конец $P_0$ регулярной области $D$
в ${\Bbb R}^n$, $n\geqslant 2,$ содержит цепь разрезов $\sigma_k$,
лежащую на сферах $S_k$ с центром в некоторой точке $x_0\in\partial
D$ и с евклидовыми радиусами $r_k\rightarrow 0$ при
$k\rightarrow\infty$. Пусть $D_k$ -- области, ассоциированные с
разрезами $\sigma_k$, $k=1,2,\ldots$. Поскольку последовательности
$x_k$ и $x_k^{\,\prime}$ сходятся к простому концу $P_0$ при
$k\rightarrow\infty,$ мы можем считать, что точки $x_k,
x^{\,\prime}_k\in D_k$ при всех $k=1,2,\ldots, .$ Соединим точки
$x_k$ и $x_k^{\,\prime}$ кривой $\gamma_k,$ полностью лежащей в
$D_k.$ Можно также считать, что континуум $A,$ относящийся к
определению класса $\frak{F}_{Q, A, p, \delta}(D),$ не пересекается
ни с одной из областей $D_k$ и что ${\rm dist}\,(\partial D,
A)>\varepsilon_0.$

\medskip
Обозначим через $C_k$ образ кривой $\gamma_k$ при отображении $f_k.$
Из соотношения (\ref{eq5E}) вытекает, что
\begin{equation}\label{eq3}
h(C_k)\geqslant a/4\qquad\forall\, k\in {\Bbb N}\,,
\end{equation}
где $h$ обозначает хордальный диаметр множества.

\medskip
Пусть $\Gamma_k$ -- семейство кривых, соединяющих $\gamma_k$ и $A$ в
$D.$ Из определения кольцевого $Q$-отображения относительно
$p$-модуля в точке $x_0$ вытекает, что
\begin{equation}\label{eq14}
M_p(f_k(\Gamma_k))\leqslant \int\limits_{A(x_0, r_k, \varepsilon_0)}
Q(x)\cdot \eta^p(|x-x_0|)\ dm(x)
\end{equation}
для каждой измеримой функции $\eta:(r_k, \varepsilon_0)\rightarrow
[0,\infty ],$ такой что $\int\limits_{r_k}^{\varepsilon_0}\eta(r)dr
\geqslant 1.$
Заметим, что функция
$$\eta(t)=\left\{
\begin{array}{rr}
\psi(t)/I(r_k, \varepsilon_0), &   t\in (r_k,
\varepsilon_0),\\
0,  &  t\in {\Bbb R}\setminus (r_k, \varepsilon_0)\,,
\end{array}
\right. $$ где $I(\varepsilon,
\varepsilon_0):=\int\limits_{\varepsilon}^{\varepsilon_0}\psi(t)dt,$
удовлетворяет условию нормировки вида (\ref{eq28*}) при $r_k$ и
$\varepsilon_0$ вместо $r_1$ и $r_2,$ поэтому из условий
(\ref{eq3.7.2}) и (\ref{eq14}) вытекает, что
\begin{equation}\label{eq11A}
M_p(f_k(\Gamma_k))\leqslant \alpha(r_k)\rightarrow 0
\end{equation}
при $k\rightarrow \infty,$ где $\alpha(\varepsilon)$ -- некоторая
неотрицательная функция, стремящаяся к нулю при
$\varepsilon\rightarrow 0,$ которая существует ввиду условия
(\ref{eq3.7.2}).

\medskip
С другой стороны, заметим, что $f_k(\Gamma_k)=\Gamma(C_k, f_k(A),
D_k^{\,\prime}),$ где $D_k^{\,\prime}=f_k(D).$ Поскольку по условию
леммы $h(f_k(A))\geqslant \delta$ при всех $k\in {\Bbb N},$ ввиду
(\ref{eq5.1}) $h(f_k(A))\geqslant \delta_1$ и
$h(C_k)\geqslant\delta_1,$ где $\delta_1:=\min\{\delta, a/4\}.$
Воспользовавшись тем, что области $D_k^{\,\prime}$ являются
равностепенно равномерными относительно $p$-модуля, мы заключаем,
что существует $\sigma>0$ такое, что
$$M_p(f_k(\Gamma_k))=M_p(\Gamma(C_k, f_k(A),
D_k^{\,\prime}))\geqslant \sigma\qquad\forall\,\, k\in {\Bbb N}\,,$$
что противоречит условию (\ref{eq11A}). Полученное противоречие
указывает на то, что предположение об отсутствии равностепенной
непрерывности семейства $\frak{F}_{Q, A, p, \delta}(\overline{D})$
было неверным. Полученное противоречие завершает доказательство
леммы.~$\Box$
\end{proof}

\medskip
\begin{lemma}\label{lem4} {\sl\, Предположим, $p\in (n-1, n],$
область $D$ регулярна и что области $D_f^{\,\prime}=f(D)$ являются
ограниченными равностепенно равномерными относительно $p$-модуля по
$\frak{R}_{Q, \delta, p, E}(D)$ областями с локально квазиконформной
границей. Пусть при $p=n$ множество $E$ имеет положительную ёмкость,
а при $n-1<p<n$ является произвольным замкнутым множеством.
Предположим также, что для каждой точки $x_0\in \overline{D}$
найдётся $\varepsilon_0=\varepsilon_0(x_0)>0$ и измеримая по Лебегу
функция $\psi(t):(0, \varepsilon_0)\rightarrow [0,\infty]$ со
следующим свойством: для любого $\varepsilon\in(0, \varepsilon_0)$
выполнено условие (\ref{eq7***}) и, кроме того, при
$\varepsilon\rightarrow 0$ выполнено условие (\ref{eq3.7.2}). Тогда
каждое из отображений $\frak{R}_{Q, \delta, p, E}(D)$ имеет
непрерывное продолжение в $\overline{D}$ и семейство $\frak{R}_{Q,
\delta, p, E}(\overline{D}),$ состоящее из всех, таким образом,
продолженных отображений $\overline{f}: \overline{D}_P\rightarrow
\overline{{\Bbb R}^n},$ является равностепенно непрерывным в
$\overline{D}_P.$ }
\end{lemma}

\medskip
\begin{proof}
Равностепенная непрерывность внутри области $D$ вытекает из
\cite[лемма~3.6.1]{Sev$_7$} в случае $p=n$ и \cite[лемма~2.4]{GSS}
при $n-1<p<n,$ а возможность продолжения каждого элемента $f$
семейства отображений $\frak{R}_{Q, \delta, p, E}(D)$ до
непрерывного отображения в замыкании $D_P$ ~--- из
\cite[лемма~3]{Sev$_5$}.

Осталось показать, что семейство $\frak{R}_{Q, \delta, p,
E}(\overline{D})$ равностепенно не\-прерывно в точках
$\partial_PD:=\overline{D}_P\setminus D.$ Предположим противное.
Рассуждая также, как и при доказательстве леммы \ref{lem3}, мы
построим две последовательности $x_k$ и $x_k^{\,\prime}\in D,$
сходящиеся к простому концу $P_0$ при $k\rightarrow\infty,$ для
которых верно соотношение вида (\ref{eq5E}). Соединим точки $x_k$ и
$x^{\,\prime}_k$ кривой $\gamma_k:[0,1]\rightarrow {\Bbb R}^n$
такой, что $\gamma_k(0)=x_k,$ $\gamma_k(1)=x^{\,\prime}_k$ и
$\gamma_k\in D$ при $t\in (0,1).$ Обозначим через $C_k$ образ кривой
$\gamma_k$ при отображении $f_k.$ Из соотношения (\ref{eq5E})
вытекает, что $h(C_k)\geqslant a/4$ при всех $k=1,2,\ldots .$

В силу \cite[лемма~2]{KR} простой конец $P_0$ регулярной области $D$
в ${\Bbb R}^n$, $n\geqslant 2,$ содержит цепь разрезов $\sigma_k$,
лежащую на сферах $S_k$ с центром в некоторой точке $x_0\in\partial
D$ и с евклидовыми радиусами $r_k\rightarrow 0$ при
$k\rightarrow\infty$. Пусть $D_k$ -- области, ассоциированные с
разрезами $\sigma_k$, $k=1,2,\ldots$. Поскольку последовательности
$x_k$ и $x_k^{\,\prime}$ сходятся к простому концу $P_0$ при
$k\rightarrow\infty,$ мы можем считать, что точки $x_k,
x^{\,\prime}_k\in D_k$ при всех $k=1,2,\ldots, .$

По определению семейства отображений $\frak{R}_{Q, \delta, p, E}(D)$
для всякого $f_k$ и всякой области $D^{\,\prime}_k:=f_k(D)$ найдётся
континуум $K_k\subset D^{\,\prime}_k$ такой, что $h(K_k)\geqslant
\delta$ и $h(f^{\,-1}(K_k),
\partial D)\geqslant \delta>0.$ Поскольку по условию леммы области $D^{\,\prime}_k$ являются
равностепенно равномерными относительно $p$-модуля, то ввиду
сказанного и учитывая условие (\ref{eq5.1}) мы получим, что при всех
$k=1,2,\ldots$ и некотором $b>0$ выполняется неравенство
\begin{equation}\label{eq13A}
M_p(\Gamma(K_k, C_k, D^{\,\prime}_k))\geqslant b\,.
\end{equation}
Рассмотрим семейство $\Gamma_k,$ состоящее из всех кривых $\beta:[0,
1)\rightarrow D^{\,\prime}_k,$ где $\beta(0)\in C_k$ и
$\beta(t)\rightarrow p\in C_k$ при $t\rightarrow 1.$ Пусть
$\Gamma^*_k$ -- семейство всех полных поднятий $\alpha:[0,
1)\rightarrow D$ семейства $\Gamma_k$ при отображении $f_k$ с
началом на $\gamma_k.$ Такое семейство корректно определено ввиду
\cite[теорема~3.7]{Vu}. Ввиду замкнутости отображения $f_k$ имеем:
$\alpha(t)\rightarrow f^{\,-1}_k(K_k)$ при $t\rightarrow 1,$ где
$f^{\,-1}_k(K_k)$ -- полный прообраз континуума $K_k$ при
отображении $f_k.$

Заметим, что ввиду компактности пространства $\overline{{\Bbb R}^n}$
при каждом фиксированном $\delta>0$ множество
$C_{\delta}:=\{x\in D: h(x, \partial D)\geqslant \delta\}$
является компактом в $D$ и $f_k^{\,-1}(K_k)\subset C_{\delta}.$
Ввиду \cite[лемма~1]{Smol} множество $C_{\delta}$ можно вложить в
континуум $E_{\delta},$ лежащий в области $D,$ при этом, можно
считать, что ${\rm dist}\,(x_0, E_{\delta})\geqslant \varepsilon_0$
за счёт уменьшения $\varepsilon_0,$ если это необходимо. Тогда на
основании (\ref{eq3*!!}) вытекает, что
\begin{equation}\label{eq10B}
M_p(f_k(\Gamma_k^*))\leqslant M_p(f_k(\Gamma(\gamma_k, E_{\delta},
D)))\leqslant \int\limits_{A(x_0, r_k, \varepsilon_0)} Q(x)\cdot
\eta^p(|x-x_0|)\ dm(x)
\end{equation}
для каждой измеримой функции $\eta: (r_k, \varepsilon_0)\rightarrow
[0,\infty ],$ такой что $\int\limits_{r_k}^{\varepsilon_0}\eta(r)dr
\geqslant 1.$
Заметим, что функция
$$\eta(t)=\left\{
\begin{array}{rr}
\psi(t)/I(r_k, \varepsilon_0), &   t\in (r_k,
\varepsilon_0),\\
0,  &  t\in {\Bbb R}\setminus (r_k, \varepsilon_0)\,,
\end{array}
\right. $$ где $I(\varepsilon,
\varepsilon_0):=\int\limits_{\varepsilon}^{\varepsilon_0}\psi(t)dt,$
удовлетворяет условию нормировки вида (\ref{eq28*}) при $r_k$ и
$\varepsilon_0$ вместо $r_1$ и $r_2,$ поэтому из условий
(\ref{eq3.7.2}) и (\ref{eq10B}) вытекает, что
\begin{equation}\label{eq11C}
M_p(f_k(\Gamma_k))\leqslant \alpha(r_k)\rightarrow 0
\end{equation}
при $k\rightarrow \infty,$ где $\alpha(\varepsilon)$ -- некоторая
неотрицательная функция, стремящаяся к нулю при
$\varepsilon\rightarrow 0,$ которая существует ввиду условия
(\ref{eq3.7.2}). Заметим, кроме того, что $f_k(\Gamma_k)> \Gamma_k$
и, одновременно, $f_k(\Gamma_k)\subset\Gamma_k,$ так что ввиду
\cite[теоремы~6.2, 6.4]{Va}
\begin{equation}\label{eq12A}
M_p(f_k(\Gamma_k))=M_p(\Gamma(K_k, C_k, D^{\,\prime}_k))\,.
\end{equation}
Однако, соотношения (\ref{eq11C}) и (\ref{eq12A}) в совокупности
противоречат (\ref{eq13A}). Полученное противоречие указывает на то,
что исходное предположение (\ref{eq6***}) было неверным, и, значит,
семейство отображений $\frak{R}_{Q, \delta, p, E}(\overline{D})$
равностепенно непрерывно в каждой точке $x_0\in E_D.$~$\Box$
\end{proof}

{\bf 3. Доказательство основных результатов.} Утверждение теорем
\ref{th3}, \ref{th4}, \ref{th1} и \ref{th2} непосредственно вытекает
из доказанных выше лемм \ref{lem1}--\ref{lem4} и \cite[лемма~3.1 и
детали доказательства теоремы~4.2]{GSS$_1$} (см. также
\cite[лемма~2.3.1]{Sev$_7$}).~$\Box$

КОНТАКТНАЯ ИНФОРМАЦИЯ

\medskip
\noindent{{\bf Евгений Александрович Севостьянов} \\
\noindent{{\bf Сергей Александрович Скворцов} \\
Житомирский государственный университет им.\ И.~Франко\\
кафедра математического анализа, ул. Большая Бердичевская, 40 \\
г.~Житомир, Украина, 10 008 \\ тел. +38 066 959 50 34 (моб.),
e-mail: esevostyanov2009@mail.ru}

\end{document}